\magnification=\magstep1



\catcode`\X=12\catcode`\@=11

\def\n@wcount{\alloc@0\count\countdef\insc@unt}
\def\n@wwrite{\alloc@7\write\chardef\sixt@@n}
\def\n@wread{\alloc@6\read\chardef\sixt@@n}
\def\r@s@t{\relax}\def\v@idline{\par}\def\@mputate#1/{#1}
\def\l@c@l#1X{\firstpart.#1}\def\gl@b@l#1X{#1}\def\t@d@l#1X{{}}

\def\crossrefs#1{\ifx\all#1\let\tr@ce=\all\else\def\tr@ce{#1,}\fi
   \n@wwrite\cit@tionsout\openout\cit@tionsout=\jobname.cit 
   \write\cit@tionsout{\tr@ce}\expandafter\setfl@gs\tr@ce,}
\def\setfl@gs#1,{\def\@{#1}\ifx\@\empty\let\next=\relax
   \else\let\next=\setfl@gs\expandafter\xdef
   \csname#1tr@cetrue\endcsname{}\fi\next}
\def\m@ketag#1#2{\expandafter\n@wcount\csname#2tagno\endcsname
     \csname#2tagno\endcsname=0\let\tail=\all\xdef\all{\tail#2,}
   \ifx#1\l@c@l\let\tail=\r@s@t\xdef\r@s@t{\csname#2tagno\endcsname=0\tail}\fi
   \expandafter\gdef\csname#2cite\endcsname##1{\expandafter
     \ifx\csname#2tag##1\endcsname\relax?\else\csname#2tag##1\endcsname\fi
     \expandafter\ifx\csname#2tr@cetrue\endcsname\relax\else
     \write\cit@tionsout{#2tag ##1 cited on page \folio.}\fi}
   \expandafter\gdef\csname#2page\endcsname##1{\expandafter
     \ifx\csname#2page##1\endcsname\relax?\else\csname#2page##1\endcsname\fi
     \expandafter\ifx\csname#2tr@cetrue\endcsname\relax\else
     \write\cit@tionsout{#2tag ##1 cited on page \folio.}\fi}
   \expandafter\gdef\csname#2tag\endcsname##1{\expandafter
      \ifx\csname#2check##1\endcsname\relax
      \expandafter\xdef\csname#2check##1\endcsname{}%
      \else\immediate\write16{Warning: #2tag ##1 used more than once.}\fi
      \multit@g{#1}{#2}##1/X%
      \write\t@gsout{#2tag ##1 assigned number \csname#2tag##1\endcsname\space
      on page \number\count0.}%
   \csname#2tag##1\endcsname}}
\def\multit@g#1#2#3/#4X{\def\t@mp{#4}\ifx\t@mp\empty%
      \global\advance\csname#2tagno\endcsname by 1 
      \expandafter\xdef\csname#2tag#3\endcsname
      {#1\number\csname#2tagno\endcsnameX}%
   \else\expandafter\ifx\csname#2last#3\endcsname\relax
      \expandafter\n@wcount\csname#2last#3\endcsname
      \global\advance\csname#2tagno\endcsname by 1 
      \expandafter\xdef\csname#2tag#3\endcsname
      {#1\number\csname#2tagno\endcsnameX}%
      \write\t@gsout{#2tag #3 assigned number \csname#2tag#3\endcsname\space
      on page \number\count0.}\fi
   \global\advance\csname#2last#3\endcsname by 1
   \def\t@mp{\expandafter\xdef\csname#2tag#3/}%
   \expandafter\t@mp\@mputate#4\endcsname
   {\csname#2tag#3\endcsname\lastpart{\csname#2last#3\endcsname}}\fi}
\def\t@gs#1{\def\all{}\m@ketag#1e\m@ketag#1s\m@ketag\t@d@l p
   \m@ketag\gl@b@l r \n@wread\t@gsin
   \openin\t@gsin=\jobname.tgs \re@der \closein\t@gsin
   \n@wwrite\t@gsout\openout\t@gsout=\jobname.tgs }
\outer\def\localtags{\t@gs\l@c@l}
\outer\def\globaltags{\t@gs\gl@b@l}
\outer\def\newlocaltag#1{\m@ketag\l@c@l{#1}}
\outer\def\newglobaltag#1{\m@ketag\gl@b@l{#1}}

\newif\ifpr@ 
\def\m@kecs #1tag #2 assigned number #3 on page #4.%
   {\expandafter\gdef\csname#1tag#2\endcsname{#3}
   \expandafter\gdef\csname#1page#2\endcsname{#4}
   \ifpr@\expandafter\xdef\csname#1check#2\endcsname{}\fi}
\def\re@der{\ifeof\t@gsin\let\next=\relax\else
   \read\t@gsin to\t@gline\ifx\t@gline\v@idline\else
   \expandafter\m@kecs \t@gline\fi\let \next=\re@der\fi\next}
\def\pretags#1{\pr@true\pret@gs#1,,}
\def\pret@gs#1,{\def\@{#1}\ifx\@\empty\let\n@xtfile=\relax
   \else\let\n@xtfile=\pret@gs \openin\t@gsin=#1.tgs \message{#1} \re@der 
   \closein\t@gsin\fi \n@xtfile}

\newcount\sectno\sectno=0\newcount\subsectno\subsectno=0
\newif\ifultr@local \def\ultralocal{\ultr@localtrue}
\def\firstpart{\number\sectno}
\def\lastpart#1{\ifcase#1 \or a\or b\or c\or d\or e\or f\or g\or h\or 
   i\or k\or l\or m\or n\or o\or p\or q\or r\or s\or t\or u\or v\or w\or 
   x\or y\or z \fi}

\def\resetall{\global\advance\sectno by 1\subsectno=0
   \gdef\firstpart{\number\sectno}\r@s@t}
\def\resetsub{\global\advance\subsectno by 1
   \gdef\firstpart{\number\sectno.\number\subsectno}\r@s@t}
\def\newsection#1\par{\resetall\vskip0pt plus.3\vsize\penalty-250
   \vskip0pt plus-.3\vsize\bigskip\bigskip
   \message{#1}\leftline{\bf#1}\nobreak\bigskip}
\def\subsection#1\par{\ifultr@local\resetsub\fi
   \vskip0pt plus.2\vsize\penalty-250\vskip0pt plus-.2\vsize
   \bigskip\smallskip\message{#1}\leftline{\bf#1}\nobreak\medskip}

\def\t@gsoff#1,{\def\@{#1}\ifx\@\empty\let\next=\relax\else\let\next=\t@gsoff
   \def\@@{p}\ifx\@\@@\else
   \expandafter\gdef\csname#1cite\endcsname{\relax}
   \expandafter\gdef\csname#1page\endcsname##1{?}
   \expandafter\gdef\csname#1tag\endcsname{\relax}\fi\fi\next}
\def\verbatimtags{\ifx\all\relax\else\expandafter\t@gsoff\all,\fi}

\def\(#1){\edef\dot@g{\ifmmode\ifinner(\hbox{\rm \noexpand\etag{#1}})
   \else\noexpand\eqno(\hbox{\noexpand\etag{#1}})\fi
   \else(\noexpand\ecite{#1})\fi}\dot@g}

\newif\ifbr@ck
\def\eat#1{}
\def\[#1]{\br@cktrue[\br@cket#1'X]}
\def\br@cket#1'#2X{\def\temp{#2}\ifx\temp\empty\let\next\eat
   \else\let\next\br@cket\fi
   \ifbr@ck\br@ckfalse\br@ck@t#1,X\else\br@cktrue#1\fi\next#2X}
\def\br@ck@t#1,#2X{\def\temp{#2}\ifx\temp\empty\let\neext\eat
   \else\let\neext\br@ck@t\def\temp{,}\fi
   \def\teemp{#1}\ifx\teemp\empty\else\rcite{#1}\fi\temp\neext#2X}
\def\resetbr@cket{\gdef\[##1]{[\rtag{##1}]}}
\def\references{\resetbr@cket\newsection References\par}

\newtoks\symb@ls\newtoks\s@mb@ls\newtoks\p@gelist\n@wcount\ftn@mber
    \ftn@mber=1\newif\ifftn@mbers\ftn@mbersfalse\newif\ifbyp@ge\byp@gefalse
\def\defm@rk{\ifftn@mbers\n@mberm@rk\else\symb@lm@rk\fi}
\def\n@mberm@rk{\xdef\m@rk{{\the\ftn@mber}}%
    \global\advance\ftn@mber by 1 }
\def\rot@te#1{\let\temp=#1\global#1=\expandafter\r@t@te\the\temp,X}
\def\r@t@te#1,#2X{{#2#1}\xdef\m@rk{{#1}}}
\def\b@@st#1{{$^{#1}$}}\def\str@p#1{#1}
\def\symb@lm@rk{\ifbyp@ge\rot@te\p@gelist\ifnum\expandafter\str@p\m@rk=1 
    \s@mb@ls=\symb@ls\fi\write\f@nsout{\number\count0}\fi \rot@te\s@mb@ls}
\def\byp@ge{\byp@getrue\n@wwrite\f@nsin\openin\f@nsin=\jobname.fns 
    \n@wcount\currentp@ge\currentp@ge=0\p@gelist={0}
    \re@dfns\closein\f@nsin\rot@te\p@gelist
    \n@wread\f@nsout\openout\f@nsout=\jobname.fns }
\def\m@kelist#1X#2{{#1,#2}}
\def\re@dfns{\ifeof\f@nsin\let\next=\relax\else\read\f@nsin to \f@nline
    \ifx\f@nline\v@idline\else\let\t@mplist=\p@gelist
    \ifnum\currentp@ge=\f@nline
    \global\p@gelist=\expandafter\m@kelist\the\t@mplistX0
    \else\currentp@ge=\f@nline
    \global\p@gelist=\expandafter\m@kelist\the\t@mplistX1\fi\fi
    \let\next=\re@dfns\fi\next}
\def\symbols#1{\symb@ls={#1}\s@mb@ls=\symb@ls} 
\def\bigsymbol{\textstyle}
\symbols{\bigsymbol\ast,\dagger,\ddagger,\sharp,\flat,\natural,\star}
\def\ftnumbers{\ftn@mberstrue} \def\ftsymbols{\ftn@mbersfalse}
\def\paginal{\byp@ge} \def\resetftnumbers{\ftn@mber=1}
\def\ftnote#1{\defm@rk\expandafter\expandafter\expandafter\footnote
    \expandafter\b@@st\m@rk{#1}}

\long\def\jump#1\endjump{}
\def\ssum{\mathop{\lower .1em\hbox{$\textstyle\Sigma$}}\nolimits}

\def\qed{\nobreak\kern 1em \vrule height .5em width .5em depth 0em}
\def\newneq{\hbox{\rlap{\hbox to 1\wd9{\hss$=$\hss}}\raise .1em 
   \hbox to 1\wd9{\hss$\scriptscriptstyle/$\hss}}}
\def\subsetne{\setbox9 = \hbox{$\subset$}\mathrel{\hbox{\rlap
   {\lower .4em \newneq}\raise .13em \hbox{$\subset$}}}}
\def\supsetne{\setbox9 = \hbox{$\subset$}\mathrel{\hbox{\rlap
   {\lower .4em \newneq}\raise .13em \hbox{$\supset$}}}}

\def\displaylinesno#1{\displ@y \tabskip\centering
  \halign to\displaywidth{\hfil$\@lign\displaystyle{##}$\hfil
    &\llap{$\@lign##$}\tabskip\z@skip\crcr
    #1\crcr}}

\def\vbar{\mathchoice{\vrule height6.3ptdepth-.5ptwidth.8pt\kern-.8pt}
   {\vrule height6.3ptdepth-.5ptwidth.8pt\kern-.8pt}
   {\vrule height4.1ptdepth-.35ptwidth.6pt\kern-.6pt}
   {\vrule height3.1ptdepth-.25ptwidth.5pt\kern-.5pt}}
\def\f@dge{\mathchoice{}{}{\mkern.5mu}{\mkern.8mu}}
\def\b@c#1#2{{\rm \mkern#2mu\vbar\mkern-#2mu#1}}
\def\b@b#1{{\rm I\mkern-3.5mu #1}}
\def\b@a#1#2{{\rm #1\mkern-#2mu\f@dge #1}}
\def\bb#1{{\count4=`#1 \advance\count4by-64 \ifcase\count4\or\b@a A{11.5}\or
   \b@b B\or\b@c C{5}\or\b@b D\or\b@b E\or\b@b F \or\b@c G{5}\or\b@b H\or
   \b@b I\or\b@c J{3}\or\b@b K\or\b@b L \or\b@b M\or\b@b N\or\b@c O{5} \or
   \b@b P\or\b@c Q{5}\or\b@b R\or\b@a S{8}\or\b@a T{10.5}\or\b@c U{5}\or
   \b@a V{12}\or\b@a W{16.5}\or\b@a X{11}\or\b@a Y{11.7}\or\b@a Z{7.5}\fi}}

\catcode`\X=11 \catcode`\@=12

\globaltags
\input amssym.def
\input amssym.tex

\def\newsection#1\par{\resetall\vskip0pt plus.1\vsize\penalty-250
   \vskip0pt plus-.1\vsize\bigskip\bigskip
   \message{#1}\leftline{\bf#1}\nobreak\bigskip}

\def\bbr{\Bbb R}

\def\b{{\bf b}}\def\a{{\bf a}}

\def\theorem#1{\medskip\noindent {\bf Theorem \stag{#1}:} }
\def\proposition#1{\medskip\noindent {\bf Proposition \stag{#1}:} }
\def\th#1{Theorem~\scite{#1}}
\def\lemma#1{\medskip\noindent {\bf Lemma \stag{#1}:} }
\def\lem#1{Lemma~\scite{#1}}
\def\remark#1{\medskip\noindent {\bf Remark \stag{#1}:} }
\def\rem#1{Remark~\scite{#1}}

\def\proof{\medskip\noindent {\bf Proof:} }
\def\proofof#1{\medskip\noindent {\bf Proof of #1:} }

\def\vspace{height5pt&\omit&height5pt&\omit&height5pt&\omit&
  height5pt&\omit&height5pt\cr}

\def\Tr{\mathop{\rm Tr}}
\def\bbc{\Bbb C}
\def\st{\tilde s}
\def\St{\tilde E}
\def\Dt{\tilde D}
\def\Nt{\tilde N}
\def\X{\cal X}

\centerline{\bf On  D. H\"agele's approach to the Bessis-Moussa-Villani 
    conjecture}
 \medskip
\centerline{Peter S. Landweber%
    \footnote*{$\scriptstyle\rm landwebe@math.rutgers.edu$}
    and Eugene R. Speer%
    \footnote{**}{$\scriptstyle\rm speer@math.rutgers.edu$}}
 \smallskip
 \centerline{Department of Mathematics}
 \centerline{Rutgers University}
 \centerline{New Brunswick, New Jersey 08903 USA}

 \bigskip
{\parindent 22pt \narrower \noindent {\bf Abstract.} The reformulation of the
Bessis-Moussa-Villani conjecture given by Lieb and Seiringer asserts that
the coefficient $\alpha_{p,r}(A,B)$ of $t^r$ in the polynomial
\hbox{$\Tr(A+tB)^p$,} with $A,B$ positive semidefinite matrices, is
nonnegative for all $p,r$.  We propose a natural extension of a method of
attack on this problem due to H\"agele, and investigate for what values of
$p,r$ the method is successful, obtaining a complete determination when
either $p$ or $r$ is odd.


   \par}

\newsection \the\sectno. Introduction

In \[H], Daniel H\"agele gives an ingenious and simple proof that if $A$ and
$B$ are $n\times n$ positive semidefinite matrices then for $p=7$ all
coefficients $\alpha_{p,r}(A,B)$ of $t$ in the polynomial 
 $$\Tr(A+tB)^p\equiv \sum_{r=0}^p\alpha_{p,r}(A,B)t^r, \(LS)$$
 where $\Tr M$ denotes the trace of the matrix $M$, are nonnegative.  If
this result could be proved for general $p$ it would imply \[LS] a
conjecture of Bessis, Moussa, and Villani \[BMV].  On the other hand, it
was also shown in \[H] that the same method does not suffice to prove the
positivity of $\alpha_{6,3}$ (we will occasionally abbreviate
``$\alpha_{p,r}(A,B)\ge0$ for all positive semidefinite $A,B$'' as
``$\alpha_{p,r}\ge0$'').  Thus it is of interest to investigate for what
values of $p$ and $r$ the method does or does not succeed in establishing
$\alpha_{p,r}\ge0$.

In this note we give several results, both negative and positive, in this
direction. We must to some extent consider separately two possible cases,
according to the parity of $p$ and $r$, and in each of these cases we define
two related integers $k$ and $q$:
  \smallskip
 {\bf Case 1:} $p$ and $r$ are odd. Then  $ p=2k+1$, $r=2q+1$;
 \smallskip
  {\bf Case 2:} $p$ is even  and $r$ is odd. Then $ p=2k+2$, $r=2q+1$.
 \smallskip\noindent
 One further case, 
  \smallskip
 {\bf Case 3:} $p$ is odd and  $r$ even,
 \smallskip\noindent
  is included implicitly; it is easy to verify that all our results for
Case~1 imply corresponding results for Case~3, obtained by replacing $r$
with $p-r$.  We will not consider in detail the case in which both $p$ and
$r$ are even; results in this case have been obtained by Klep and
Schweighofer \[K,KS] and by Burgdorf \[B], as we discuss briefly in
Section~4.  In each of Cases~1 and 2 we define precisely a proof strategy
which is the natural generalization of that of \[H] and investigate its
success.  We are able to classify completely the pairs $(p,r)$ for which
the method succeeds; unfortunately, although these include one infinite
class ($p,r$ odd with $r=p-4$), the method does not succeed in enough cases
to establish the BMV conjecture.

Results of this sort should be viewed in the light of an important theorem
of Hillar \[Hil], which implies that if $\alpha_{p,r}\ge0$ then also
$\alpha_{p',r'}\ge0$ if $p\ge p'$, $r\ge r'$, and $p-r\ge p'-r'$.  For
example, it is pointed out in \[H] that although the proof method used
there does not apply directly when $p=6$, $r=3$, the nonnegativity of
$\alpha_{6,3}$ follows from the corresponding result for $p=7$, $r=3$;
similarly, our result that $\alpha_{p,p-4}\ge0$ for $p$ odd implies the
positivity of $\alpha_{p,r}$, for all $p$, when $r\le4$ or $r\ge p-4$.
Moreover, it follows that to establish the full BMV conjecture it suffices
to establish positivity of $\alpha_{p_n,r_n}$ for some sequences $p_n,r_n$
with $p_n\to\infty$, $r_n\to\infty$, and $p_n-r_n\to\infty$ as
$n\to\infty$.  Our results leave open the possibility of proving the BMV
conjecture by successfully applying the method of \[H] to such a sequence
with $p_n,r_n$ even.  

In order to describe the method more precisely we write $X_0\equiv A$ and
$X_1\equiv B$.  Let $E_{p,r}$ be the set of binary strings of length $p$,
$s=s_1\cdots s_p$, containing exactly $r$ 1's, and for $s\in E_{p,r}$ write
$Y_s=X_{s_1}\cdots X_{s_p}$.  Then
 $$\alpha_{p,r}(A,B)=\sum_{s\in E_{p,r}}\Tr(Y_s).\(alpha)$$
 Now for
coefficients $c=(c_u)_{u\in E_{k,q}}\in\bbc^{E_{k,q}}$ define
$Z(c)=\sum_{u\in E_{k,q}}c_uY_u$.  Then  we will have
$\alpha_{p,r}(A,B)\ge0$ if we can show that for some appropriately chosen
$c^{(m)}=(c^{(m)}_u)_{u\in E_{k,q}}$, $1\le m\le M$,
 $$\alpha_{p,r}(A,B)=\cases{
   \sum_m \Tr\bigl(Z(c^{(m)})BZ(c^{(m)})^*\bigr),& in case 1,\cr
   \sum_m \Tr\bigl(Z(c^{(m)})BZ(c^{(m)})^*A\bigr),& in case 2.\cr}
   \(alphapos)$$
This follows from the fact that if $a$ and $b$ are the nonnegative square
roots of $A$ and $B$, respectively, then
$\Tr(Z(c)BZ(c)^*)=\Tr[(Z(c)b)(Z(c)b)^*]$ and
$\Tr(Z(c)BZ(c)^*A)=\Tr[(aZ(c)b)(aZ(c)b)^*]$.

To relate \(alpha) with \(alphapos) we must make explicit the effect of the
invariance of the trace under cyclic permutations.  Let $\St_{p,r}$ be the
set of equivalence classes of $E_{p,r}$ modulo cyclic permutations, with
$\pi:E_{p,r}\to \St_{p,r}$ the canonical projection.  Then \(alpha) becomes
 $$\alpha_{p,r}(A,B) 
   =\sum_{\st\in \St_{p,r}}|\st| \Tr(Y_{s(\st)}),\(newalpha)$$
 where $|\st|$ is the number of elements in $\st$ and $s(\st)$ is some
element of $\st$.  Similarly, if we define
$\sigma\equiv\sigma_{p,r}:E_{k,q}\to E_{p,r}$ by
 $$\sigma_{p,r}(u,v)=\cases{u_1\cdots u_k1v_k\cdots v_1,&in case 1,\cr
u_1\cdots u_k1v_k\cdots v_10,&in case 2,\cr}
$$
  then the right hand side of \(alphapos) becomes
 $$\sum_{\st\in \St_{p,r}}\sum_{(u,v)\in(\pi\sigma)^{-1}(\st)}\sum_m
            c^{(m)}_u\overline c^{(m)}_v \Tr(Y_{\sigma(u,v)}),\(rhs)$$
 so that \(alphapos) will hold for all $A,B$ if for all $\st\in\St_{p,r}$,
 $$|\st|=\sum_{(u,v)\in(\pi\sigma)^{-1}(\st)}\sum_m
            c^{(m)}_u\overline c^{(m)}_v.\(key)$$
 The generalization of the method of \[H] referred to above is establish
 
 \smallskip\noindent
 {\bf Condition H:} There exist $M\ge1$ and coefficients $c^{(m)}$,
$m=1,\ldots,M$, such that \(key) is satisfied for all $\st\in\St_{p,r}$.

 \smallskip
 Before proceeding we verify a fact which is obviously necessary for the
existence of such $c^{(m)}$.

\proposition{nec} {\sl For any $p$ and $r$ and any $\st\in\St_{p,r}$
there exist $u,v\in \St_{k,q}$ such that $\sigma_{p,r}(u,v)\in\st$.}

\proof We give the proof in Case 1;  Case~2 is similar.  
A useful geometric picture (the reader might draw a sketch) is
obtained by letting $C\subset\bbc$ denote the set of $p^{\rm th}$ roots of
unity and identifying an element $s=s_1\cdots s_p\in E_{p,r}$ with a map
$s:C\to\{0,1\}$ labeling the elements of $C$; the identification is via
$s(\exp 2j\pi i/p)=s_j$, $j=1,\ldots,p$.  Any $\theta\in\bbr$ defines the
line $L_\theta$ in $\bbc$ through the origin and the point
$z_\theta=\exp i\theta$, oriented from the origin toward $z_\theta$.  Let
$N_1(\theta)$ be the number of points $\omega\in C$ for which $s(\omega)=1$
and which lie to the right of $L_\theta$, let $N_2(\theta)$ be the number
of such points which lie to the left of $L_\theta$, and let
$N(\theta)=N_1(\theta)-N_2(\theta)$.  $N(\theta)$ is odd unless
$\pm z_\theta\in C$ with $s(\pm z_\theta)=1$, in which case it is even, and
if $N(\theta_0)=0$ for some $\theta_0$ then we can immediately read off the
desired $u,v$.  But taking $\theta$ with $\pm z_\theta\notin C$, so that
$N(\theta)$ is odd, we observe that $N(\theta+\pi)=-N(\theta)$ and so
$N(\theta_0)=0$ for some intermediate $\theta_0$.\qed

\newsection \the\sectno. Positive results

In this section we show that Condition~H holds in the following cases:
 \smallskip\noindent
 {\bf Case 1:} $r=1$; $r=p-2$; $r=p-4$; and $p=11,r=3$.  The cases $r=1$
and $r=p-2$ are easy (in each case one takes $M=1$ and $c^{(1)}_u=1$ for all
$u\in E_{k,q}$); the remaining cases are
covered in Theorems~\scite{p11r3} and \scite{pminus4} below.
 \smallskip\noindent
 {\bf Case 2:} $r=1$ and $r=p-1$.  These follow the pattern
 of the two easy cases above; verification is left to the reader. 

 \theorem{p11r3} {\sl Condition~H holds if $p=11$ and $r=3$.}

\proof  Defining
 $$\eqalign{Z_1&=Y_{00001}+Y_{00010}+Y_{00100}+Y_{01000}-Y_{10000},\cr
Z_2&=\sqrt2\,(Y_{00100}-Y_{01000}-Y_{10000}),\cr
Z_3&=2\,(Y_{00100}-Y_{01000}),\cr
Z_4&=2\,Y_{01000},\cr}\(51coeff)$$
 and using the fact that, since $p=11$ is prime, $|\st|=11$ for all $\st\in
 \St_{11,3}$, one finds easily that (compare \(alphapos), case 1) 
 $$\alpha_{11,3}(A,B) =11 \sum_{i=1}^4 \Tr( Z_iBZ_i^*).\qed$$

   We remark that both positive and negative coefficients occur among the
$c^{(m)}_u$ implicitly defined by \(51coeff). It can easily be shown that
no solution in which all the coefficients are positive is possible; this is
in contrast to the situation for the case $p=7$, $r=3$ discussed in \[H]
and for the cases treated in  \th{pminus4} below.

\theorem{pminus4} {\sl Condition~H holds if $p$ is odd and $r=p-4$.}
 
 \medskip
  Note that the case $p=7$, $r=3$ of this theorem appears in \[H]; the case
$p=9$, $r=5$ was obtained by Klep and Schweighofer (see \[K]).  After we had
completed our work we learned that \th{pminus4} was obtained independently
by Burgdorf \[B]. 

 The theorem will follow almost immediately from the next lemma.

\lemma{partition} {\sl Let $p=2k+1\ge5$ and let $r=p-4=2q+1$.  Then
$E_{k,q}$ may be partitioned as $E_{k,q}=\bigcup_{m=1}^{k-1}D_m$ in such a
way that for every $\st\in\St_{p,r}$ there exists a unique $m$, $1\le m\le
k-1$, and unique $u,v\in D_m$, such that $\sigma(u,v)\in\st$.}

\proofof{\th{pminus4}} Set $p=2k+1$ and $p-4=2q+1$.  We must find
coefficients $c^{(m)}=(c^{(m)}_u)_{u\in E_{k,q}}$ satisfying \(key); since
$p$ and $p-4$ are relatively prime, $|\st|=p$ for every $\st\in\St_{p,r}$
and so equivalently we must find $c^{(m)}$ satisfying
 $$\sum_m\sum_{(u,v)\in(\pi\sigma)^{-1}(\st)}
            c^{(m)}_u\overline c^{(m)}_v = 1.\(key1)$$
 But from \lem{partition}, \(key1) holds if $c^{(m)}$, $m=1,\ldots,k-1$,
is the characteristic function of $D_m$: $c^{(m)}_u=1$ if $u\in D_m$,
$c^{(m)}_u=0$ otherwise.\qed

 \medskip
 The next proof is somewhat complicated;  it might help the
reader to work through it in the case $p=9$, $r=5$ (this was the case that
suggested the general result).

\proofof{\lem{partition}} Recalling that an element $u\in E_{k,q}$ is a binary
string $u_1\,u_2\,\cdots\,u_k$, we define
 $$\eqalign{D_1&=\{u\in E_{k,q}\mid u_1=0\},\cr
D_2&=\{u\in E_{k,q}\mid u_1=1, u_k=0\},\cr
D_3&=\{u\in E_{k,q}\mid u_1=u_k=1, u_2=0\},\cr
D_4&=\{u\in E_{k,q}\mid u_1=u_k=u_2=1,u_{k-1}=0\},\qquad \rm etc.,\cr}
$$
 and in general, for $j\ge0$,
 $$\displaylines{D_{2j+1}=\{u\in E_{k,q}\mid 
    u_1=u_2=\cdots=u_j=u_k=u_{k-1}=\cdots=u_{k-j+1}=1,\ u_{j+1}=0\},\cr
 D_{2j+2}=\{u\in E_{k,q}\mid 
    u_1=u_2=\cdots=u_{j+1}=u_k=u_{k-1}=\cdots=u_{k-j+1}=1,\ u_{k-j}=0\}.\cr}
$$
 It is clear that the $D_m$ so defined form a partition of $E_{k,q}$. We
 will write $\Dt_m=\sigma(D_m\times D_m)$, so that we must prove that for
 any $\st\in \St_{p,r}$, $|\st\cap\bigcup_{m=1}^{k-1}\Dt_m|=1$. 

Note that a string $u\in E_{k,q}$ contains exactly two zeros, and if
$u\in D_m$ then the position of one of these zeros is fixed and there are
$k-m$ possible positions for the remaining one; thus $|D_m|=k-m$.  Note
also that $u,v\in D_m$ if and only if the form of $\sigma(u,v)$ is
 $$\eqalignno{1^j\,0\,w\,1^j\,1\,1^j\,x\,0\,1^j\,,&\quad
     \hbox{if $m=2j+1$, $j\ge0$,}&\(form/a) \cr
  1^{j+1}\,w\,0\,1^j\,1\,1^j\,0\,x\,1^{j+1}\,,&\quad
     \hbox{if $m=2j+2$, $j\ge0$}&\(form/b) \cr} $$
 where $w,x\in E_{k-m,k-m-1}$ are arbitrary.  

Now fix $\st\in\St_{p,r}$.  There are nonnegative integers
$n_0,\ldots,n_3$, with $n_0+n_1+n_2+n_3=2k-3$, such that $\st$ consists of
all cyclic permutations of the string
 $$0\,1^{n_0}\,0\,1^{n_1}\,0\,1^{n_2}\,0\,1^{n_3}.\(string)$$
  We must show that precisely one element of $\st$ has one of the forms
\(form).  

Consider first \(form/a); the initial $1^j\,0$ and final $0\,1^j$ there
imply that if that string is put in the form \(string) by a cyclic
permutation then it will contain a substring $0\,1^{2j}\,0$, i.e., that if
an element in $\st$ has the form \(form/a) then one of the integers $n_i$
must be even.  Conversely, if $n_i$ is even for some $i$, with $n_i=2j_i$
($j_i\ge0$), then the string $s_i\in\st$ defined by
 $$s_i=1^{j_i}\,0\,1^{n_{i+1}}\,0\,1^{n_{i+2}}\,0\,1^{n_{i+3}}\,0\,1^{j_i}
   \(cand/a)$$
 (here addition on the indices of the $n_l$'s is taken modulo 4) will lie
in $\Dt_{n_i+1}$ if $n_{i+1}$ and $n_{i-1}$ satisfy certain
additional constraints, which we discuss below. The discussion of \(form/b)
is similar: if some $n_i$ is odd, $n_i=2j_i+1$ ($j_i\ge0$), then the cyclic
permutation of \(string) in which the block $1^{n_i}$ is moved to the
center is a candidate to lie in $\Dt_{n_i+1}$. If $j_i+n_{i-1}+2\le k$ and
$j_i+n_{i+1}+2\le k$ (the only case that will be relevant, since \(form/b)
has two zeros on each side of its center) then this string
has the form
 $$s_i=  1^{k-(j_i+n_{i-1}+2)}\,0\,1^{n_{i-1}}\,0\,1^{j_i}\,1\,1^{j_i}\,0
      1^{n_{i+1}}\,0\,1^{k-(j_i+n_{i+1}+2)},
\(cand/b)$$
 and will lie in $\Dt_{n_i+1}$ under further constraints on $n_{i\pm1}$.
We see that for each $i$, $i=0,1,2,3$, there is one possible element of
$\st$ which could lie in $\Dt_{n_i+1}$, given by \(cand/a) or \(cand/b) as
$n_i$ is even or odd.

Now we ask what further conditions on $n_{i\pm1}$
would imply that \(cand/a) has the form \(form/a) or \(cand/b) the form
\(form/b). Consider first \(cand/a), and recall that here $n_i=2j_i$. The
second zero in \(cand/a) is located at position $j_i+n_{i+1}+2$, and for
\(cand/a) to have the form \(form/a) it is necessary that this zero lie to
the left of a block $1^{j_i}\,1\,1^{j_i}$ at the center of the string, that
is, to the left of position $k-j_i+1$.  Thus $s_i\in \Dt_{n_i+1}$ is possible
only if $j_i+n_{i+1}+2<k-j_i+1$, i.e., only if $n_i+n_{i+1}\le k-2$.
Combining this result with that of a similar analysis of the position of
the third zero shows that
 $$\hbox{$s_i\in \Dt_{n_i+1}$ if and only if $n_i+n_{i+1}\le k-2$
   and $n_i+n_{i-1}\le k-2$.}\(cond)$$
 The analysis of \(cand/b), where $n_i=2j_i+1$, is similar: for this to have
the form \(form/b), there must be at least $j_i+1$ initial ones in the
string, requiring that $k-(j_i+n_{i-1}+2)\ge j_i+1$; since there must also
be $j_i+1$ ones at the end of the string we are led again to the
conclusion \(cond).

Finally we observe that the condition that $\sum_{i=0}^3n_i=2k-3$ implies
that of any pair of inequalities $n_i+n_{i+1}\le k-2$ and
$n_{i+2}+n_{i+3}\le k-2$ exactly one must be true. This implies that the
condition of \(cond) will be satisfied for exactly one value of $i$ (modulo
4), so that $s_i\in \Dt_{n_i+1}$ (that is,
$\st\cap\Dt_{n_i+1}=\{s_i\}$) holds for precisely one value of $i$.  From
\(cand/a) or \(cand/b) one can then read off the unique $u,v\in D_{n_i+1}$
such that $\sigma(u,v)=s_i$.
\qed

 \newsection \the\sectno. Negative results

In this section we show that Condition~H does not hold in the following cases:
 \smallskip\noindent
 {\bf Case 1:} $5\le r\le p-6$; $p\ge13$, $r=3$; and $p=9$, $r=3$.
 \smallskip\noindent
 {\bf Case 2:} $3\le r\le p-3$.
 \smallskip\noindent
The method of proof in all of these cases is similar to the argument of
\[H] establishing a negative result for $p=6$, $r=3$.

Throughout the rest of this section we assume that we are in case~1 or
case~2, that is, that $r=2q+1$ is odd, but to the extent possible we treat
these two cases in a unified manner, so that for the moment either $p=2k+1$
or $p=2k+2$.  If $u,v\in E_{k,q}$ we write $\Nt(u,v)=|\pi(\sigma(u,v))|$
and $N(u,v)=|(\pi\sigma)^{-1}(\pi(\sigma(u,v)))|$; that is, $\Nt(u,v)$ is
the number of distinct strings obtained from $\sigma(u,v)$ by cyclic
permutation, and $N(u,v)$ is the number of ordered pairs
$(w,x)\in E_{k,q}\times E_{k,q}$ such that $\sigma(w,x)$ is obtained from
$\sigma(u,v)$ by a cyclic permutation.  We will compute $N(u,v)$ using the
following simple remark.

\remark{size} Let $k'=p-k-1$ so that $k'=k$ in case~1, $k'=k+1$ in case~2.
Then for any $s\in E_{p,r}$ with $|\pi(s)|=p$, $|(\pi\sigma)^{-1}(\pi(s))|$
is equal to the number of indices $i$, $1\le i\le p$, such that (i)~$s_i=1$
and the preceding (if $i\ge k'+1$) or succeeding (if $i\le p-k'$) $k'$
entries of $s$---that is $s_{i-k'}\cdots s_{i-1}$ or
$s_{i+1}\cdots s_{i+k'}$, respectively---contain exactly $q$ ones, and (ii)
in case~2, if also $i-k'=0$ or $i+k'=0$, respectively.  Of course if
$s=\sigma(u,v)$ then $i=k+1$ satisfies this criterion.  The application of
this remark in any particular case is straightforward but tedious; we give
a full discussion of one case in the proof of \lem{first} and after that we
are rather sketchy, leaving the details to the reader.  It is probably most
helpful to work out a simple example in each case.

  \medskip
  We now define
$w=0^{k-q}\,1^q\in E_{k,q}$.

\lemma{first} {\sl Suppose that $u\in E_{k,q}$. Then (a) $\Nt(w,u)=p$, and
(b)~if $u_1=0$ or $p$ is even (i.e., we are in case~2) then $N(w,u)=1$.  In
particular, (c)~$\Nt(w,w)=p$ and $N(w,w)=1$.}

\proof (a) The string $\sigma(w,u)$ contains a substring of at least $q+1$
consecutive ones, and since there are a total of $2q+1$ ones in the string, no
nontrivial cyclic permutation of $\sigma(w,u)$ can coincide with it.

 \smallskip\noindent
 (b) Under either hypothesis, $s\equiv\sigma(w,u)$ has the form
$s=0^{k-q}\,1^q\,1\,s_{k+2}\,\cdots\,s_{p-1}\,0$; the key observation is
that for $1\le j\le q+1$ the last $j$ entries of $s$ can contain at most
$j-1$ ones, and so entries $k+2,\ldots,p-j$ must contain at least $q-j+1$
ones. We show that no index $i$, $1\le i\le p$, other than $i=k+1$, can
satisfy criterion (i) of Remark~\scite{size}. Suppose then that $s_i=1$ and
$i\ne k+1$.  There are three possible cases: if $k-q+1\le i\le k$ then
$s_{i+1}\cdots s_{i+k'}=1^{k-i}\,1\,s_{k+2}\,\cdots \,s_{p-(k+1-i)}$
contains, by the observation above, at least $(k-i)+1+(q-k+i)=q+1$ ones; if
$k+2\le i\le p-q-1$ then $s_{i-k'}\,\cdots\,s_{i-1}$ contains the substring
$s_{k-q+1}\,\cdots\,s_{k+1}=1^{q+1}$; and if $p-q\le i\le p$
then $s_{i-k'}\,\cdots\,s_{i-1}=1^{p-i}\,1\,s_{k+2}\,\cdots\,s_{i-1}$
contains at least $(p-i)+1+(q-p+i)=q+1$ ones.

 \smallskip\noindent
 (c) This is an immediate consequence of (a) and (b).\qed

\lemma{neg} {\sl Suppose there exist $x,y,z\in E_{k,q}$, all distinct from
$w$ and with $x\ne y$ and $x\ne z$, such that
 $$\displaylinesno{N(w,x)=N(w,y)=N(x,x)=1,&\(lneg/a)\cr
     \Nt(w,x)=\Nt(w,y)=\Nt(x,x)=p,&\(lneg/b)\cr
     N(z,z)=3,\hbox{ with }
     \pi(\sigma(z,z))=\{\sigma(z,z),\sigma(x,y),\sigma(y,x)\}.&\(lneg/c)\cr}$$
 Then Condition~H does not hold.}

 \medskip
 We remark that the requirement that all of $x$, $y$, $z$ and $w$ be
distinct, except for the possibility that $y=z$, actually follows from
\(lneg) and \lem{first}.

\proof We suppose that for some $M$ and $c^{(m)}$, \(key) holds for all
$\st$, and derive a contradiction.  From \(key) applied to
$\pi(\sigma(w,w))$, $\pi(\sigma(x,x))$, and $\pi(\sigma(w,x))$
we have, using \lem{first}(c) and \(lneg/a)--\(lneg/b),
 $$p=\sum_mc^{(m)}_{w}\overline c^{(m)}_{w}
   = \sum_mc^{(m)}_{x}\overline c^{(m)}_{x}
   = \sum_mc^{(m)}_{w}\overline c^{(m)}_{x}.\(firstii)$$
 These equations, together with the standard necessary condition for
equality to hold in the Cauchy-Schwarz inequality, then imply that
 $$c^{(m)}_{w}= c^{(m)}_{x}, \qquad \hbox{for $m=1,\ldots,M.$}\(secondii)$$
 But, first from \(key) applied to $\pi(\sigma(w,y))$, and then from
\(secondii),
 $$p=\sum_mc^{(m)}_{w}\overline c^{(m)}_{y}
   = \sum_mc^{(m)}_{x}\overline c^{(m)}_{y}.\(thirdii)$$
 Finally, from \(key) applied to $\pi(\sigma(z,z))$, \(lneg/c), and then
\(thirdii),
 $$\Nt(z,z)=\sum_mc^{(m)}_{z}\overline c^{(m)}_{z}
   + \sum_mc^{(m)}_{x}\overline c^{(m)}_{y}
   + \sum_mc^{(m)}_{y}\overline c^{(m)}_{x}
  = \sum_mc^{(m)}_{z}\overline c^{(m)}_{z}+2p\ge2p,\(finalii)$$
 a contradiction, since $\Nt(z,z)$ must divide $p$.\qed

\theorem{negative} {\sl If $r$ is odd and (a)~$p$ is odd and
$5\le r\le p-6$, (b)~$p$ is odd, $p\ge13$, and $r=3$, or (c)~$p$ is even
and $3\le r\le p-3$, then Condition~H does not hold.}

\proof  (a) In this case we claim that the strings
 $$x=0\,1\,0^{k-q-1}\,1^{q-1},\quad y=0^{k-q-2}\,1^q\,0^2,\quad
   \hbox{and} \quad z=0\,1^q\,0^{k-q-1},$$
 fulfill the conditions of \lem{neg}. Since $2\le q\le k-3$ we have
$x\ne y$ and $x\ne z$ (although $y=z$ if $q=k-3$).  The
conditions
 $$N(w,x)=N(w,y)=1, \qquad
     \Nt(w,x)=\Nt(w,y)=p,$$
 follow from \lem{first}, since $x_1=y_1=0$.  

Consider now
$\sigma(x,x)=0\,1\,0^{k-q-1}\,1^{q-1}\,1\,1^{q-1}\,0^{k-q-1}\,1\,0$;
this contains a unique string of $2q-1\ge3$ consecutive ones and so can
never coincide with a cyclic permutation of itself, so that indeed
$\Nt(x,x)=p$. A detailed analysis using \rem{size}, as in the proof of
\lem{first}(b) (but by symmetry it is necessary to consider only $i\le k$),
shows that $N(x,x)=1$.

Finally consider
$s\equiv\sigma(z,z)=0\,1^q\,0^{k-q-1}\,1\,0^{k-q-1}\,1^q\,0$.  Again,
consideration of the sizes of the three blocks of consecutive ones shows
that $\Nt(z,z)=p$.  To find $N(z,z)$ we note that a cyclic
permutation which brings the one at position $i=2$ of $s$ to the center
position $i=k+1$ yields that string
$0^{k-q-2}\,1^q\,0^2\,1\,1^{q-1}\,0^{k-q-1}\,1\,0=\sigma(y,x)$, and one
obtains $\sigma(x,y)$ by a cyclic permutation bringing the one at
$i=p-1$ in $s$ to $i=k+1$.  However, if $3\le i\le q+1$ then
$s_{i+1}\cdots s_{i+k}$ contains at most $q-1$ ones, with a similar
conclusion if $p-q-1\le i\le p-2$, so that $N(z,z)=3$ and \(lneg/c)
holds.

 \smallskip\noindent
 (b) In this case the strings
 $$x=0^{k-3}\,1\,0^2,\quad  y=z=0\,1\,0^{k-2},$$
 fulfill the conditions of \lem{neg}; the verification is similar
 to the above.

 \smallskip\noindent
 (c) If $5\le r\le p-3$ then the strings
 $$x=10^{k-q}\,1^{q-1},\quad  y=0^{k-q-1}\,1^q\,0,\quad z=1^q0^{k-q},$$
 fulfill the conditions of \lem{neg}; again the verification is similar
 to that of case (a). If $r=3$  and  $p\ge8$ then the conclusion  follows
 from the case $r=p-3$ after the interchange of $A$ and $B$.  Finally, the
 result for case $p=6$, $r=3$ was established in \[H].\qed

 \medskip
 The next result covers the one remaining negative result not included in
\th{negative}.  It is stated without proof in \[K].

\theorem{negative93} {\sl If $p=9$ and $r=3$ then Condition~H does not
hold.}

\proof Again  we suppose that there exist $c^{(m)}$, $m=1,\ldots,M$,
so that \(key) holds for all $\st$, and derive a contradiction by looking
at a few specific choices of $\st$, as  given in Table
1; there we write $v_1=0001$, $v_2=0100$ (with $v_1,v_2\in E_{4,1}$).
\bigskip

\centerline{\vbox{\hrule\offinterlineskip%
 \halign{\vrule#\tabskip10pt&\strut\hfill$#$\hfill&\vrule#&
   \hfill$#$\hfill&\vrule#&
  \hfill$#$\hfill&\vrule#&\hfill$#$\hfill&\vrule#\tabskip0pt\cr
\vspace
&\hbox{ Name of }\st&&\hbox{Typical } s\in \st      &&
          |\st|&&(\pi\sigma)^{-1}(\st)&\cr
\vspace
\noalign{\hrule}
\vspace 
&\st_1   &&000111000   &&9 &&\{(v_1,v_1)\}&\cr
\vspace 
&\st_2   &&010010010   &&3 &&\{(v_2,v_2)\}&\cr
\vspace 
&\st_3   &&000110010   &&9 &&\{(v_1,v_2)\}&\cr
\vspace 
}\hrule}}
 \medskip
\centerline{Table 1}
 \bigskip\noindent
 From \(key) applied to $\st_1$, $\st_2$, and $\st_3$, we have
 $$9=\sum_mc^{(m)}_{v_1}\overline c^{(m)}_{v_1}
   = \sum_mc^{(m)}_{v_1}\overline c^{(m)}_{v_2};\qquad
   3= \sum_mc^{(m)}_{v_2}\overline c^{(m)}_{v_2}.\(first)$$
 These equations, however, are inconsistent with the Cauchy-Schwarz
inequality.\qed

 \newsection 4. Concluding remarks

In recent work \[KS] Klep and Schweighofer give a systematic algebraic
language in which to discuss the method of \[H].  They introduce the
associative $\bbr$-algebra $\bbr\langle \a,\b\rangle$ with noncommuting
generators $\a$ and $\b$ ($X$ and $Y$ in the notation of \[KS]), furnished
with a natural involution $f\mapsto f^*$ obtained by reversing each word in
the generators. They further define
$\Sigma^2\subset\bbr\langle\a,\b\rangle$ to be the cone of elements
$f\in\bbr\langle\a,\b\rangle$ which may be written as sums of Hermitian
squares, $f=\sum_i g_i^*g_i$, and $\Theta^2$ to be the cone of elements
which are cyclically equivalent to elements of $\Sigma^2$, where two
elements $f$ and $g$ are cyclically equivalent if their difference is a sum
of commutators.  It follows that if $f(\a,\b)\in\Theta^2$ and $a,b$ are
nonnegative $n\times n$ matrices then $\Tr(f(a,b))\ge0$, so in order to
show that $\alpha_{p,r}\ge0$ it suffices to verify that
$S_{p,r}(\a^2,\b^2)\in\Theta^2$, where
$S_{p,r}(\a^2,\b^2)\in\bbr\langle\a,\b\rangle$ denotes the sum of all
possible products of $r$ factors $\b^2$ and $p-r$ factors $\a^2$.

It is immediate that if hypothesis $H$ is satisfied for some $p,r$ falling
under Case~1 or Case~2, or if $H$ is satisfied for $p,p-r$ with $p,r$
falling under Case~3, then $S_{p,r}(\a^2,\b^2)\in\Theta^2$; further, it
follows from a result of \[KS] (Proposition 2.2) that the converse also
holds.  This means that the results of Sections~2 and 3 establish, for
every $p,r$ with either $p$ or $r$ odd, whether or not
$S_{p,r}(\a^2,\b^2)\in\Theta^2$.  In particular, we can conclude that the
approach of \[H] (at least as formulated in \[KS]) when applied to such $p$
and $r$ cannot establish the BMV conjecture for any $p$ larger than 9.

Thus to make progress on the BMV conjecture using this approach one must
consider cases in which both $p$ and $r$ are even.  In this direction, Klep
and Schweighofer show \[KS] that $S_{14,4}(\a^2,\b^2)$ and
$S_{14,6}(\a^2,\b^2)$ belong to $\Theta^2$, which, together with results of
\[Hil] or by independent arguments given in \[KS], implies that the BMV
conjecture is satisfied for $p=13$ and indeed, by \[Hil], for $p\le13$.
Moreover, Burgdorf \[B] has obtained a version of \th{pminus4} strengthened
to include $p,r$ even: she shows that $S_{p,4}(\a^2,\b^2)\in\Theta^2$ (and
hence $S_{p,p-4}(\a^2,\b^2)\in\Theta^2$) for all $p\ge4$.

\newsection Acknowledgments

We thank I.~Klep, M.~Schweighofer, and S.~Burgdorf for communicating
their results to us in manuscript form, and C.~Hillar and D.~H\"agele for
helpful comments on a preliminary version of this paper.

\references
\parindent0pt

 \[BMV] D.~Bessis, P.~Moussa, and M.~Villani, Monotone converging
variational approximations to the functional integrals in quantum
statistical mechanics. {\it J.~Math.~Phys.}~{\bf 16}, 2318--2325 (1975). 

\smallskip
\[B] S.~Burgdorf, private communication.

 \smallskip
\[H] D. H\"agele, Proof of the cases $p\le7$ of the Lieb-Seiringer
formulation of the Bessis-Moussa-Villani conjecture.  {\it
J.~Stat.~Phys.}~{\bf 127}, 1167--1171 (2007).

 \smallskip
 \[Hil] C. J. Hillar, Advances on the Bessis-Moussa-Villani trace
conjecture. {\it Linear Alg. Appl.}~{\bf 426}, 130--142 (2007). 

 \smallskip
 \[K] I.~Klep (joint work with M. Schweighofer), Sums of Hermitian squares,
Connes' embedding problem and the BMV conjecture, {\it Mathematisches
Forschungsinstitut Oberwolfach}, Report No.  14/2007: Reelle Algebraische
Geometrie, pp. 27--30. March 11th--March 17th, 2007.

\smallskip
\[KS] I.~Klep and M.~Schweighofer, Sums of hermitian squares and the BMV
conjecture. arXiv:0710.1074.

 \smallskip
 \[LS] E.~H.~Lieb and R.~Seiringer, Equivalent forms of the
Bessis-Moussa-Villani conjecture.  {\it J.~Stat.~Phys.}~{\bf 115}, 185--190
(2004).

\bye